\documentclass[openbib,twoside,notitlepage]{article}  

\usepackage{amssymb}
\usepackage{latexsym}
\usepackage{exscale}
\usepackage{epsf}
\usepackage{theorem}
\usepackage[dvips]{graphicx}

\newfont{\bsym}{cmbsy10 scaled\magstep2}
\newfont{\bsymi}{cmbsy10}
\newfont{\bmath}{cmmib10 scaled\magstep2}
\newfont{\titfont}{cmbx10 scaled \magstep3}
\newfont{\eighrm}{cmr8}
\newfont{\ack}{cmr10}
\newfont{\Ack}{cmbx10 at 14pt}

\newtheorem{theorem}{Theorem}[section]
\newtheorem{lemma}{Lemma}[section]
\newtheorem{definition}{Definition}[section]
{\theorembodyfont{\rmfamily} 
\newtheorem{remark}{Remark}[section]}

\newcommand{\vsp}{\vspace{12pt}}

\newcommand{\R}{{I\!\! R}}
\newcommand{\NN}{{I\!\! N}}

\newcommand{\II}{{\cal I}}
\newcommand{\JJ}{{\cal J}}

\newcommand{\tv}{\hbox{Tot.Var.}}

\newcommand{\n}{\noindent}

\newcommand{\be}{\begin{equation}}
\newcommand{\ee}{\end{equation}}

\newcommand{\sgn}{\textrm{sgn }}

\begin{document}

\setlength{\voffset}{-0.5in}

\title{
{\huge Traffic Flow on a Road Network}
\vspace{0.1in}}

\author{
{\scshape Giuseppe Maria Coclite} 
\thanks{SISSA-ISAS, via Beirut 2-4, 34014 - Trieste, \ Italy; E-mail:
\texttt{coclite@sissa.it.}} 
\and    
        {\scshape Benedetto Piccoli } 
 \thanks{Istituto per le Applicazioni del Calcolo "M. Picone", Viale del
Policlinico 137,  00161 - Roma,\ Italy;\ E-mail:
\texttt{piccoli@iac.rm.cnr.it}}}
  \date{February 2002}

\maketitle
\thispagestyle{empty}

\vspace{-0.6cm}
\begin{abstract}
\vspace{5pt}

\noindent
This paper is concerned with a fluidodynamic model for traffic flow.
More precisely, we consider a single conservation law, deduced from
conservation of the number of cars, defined on a road network that is a
collection of roads with junctions.
The evolution problem is underdetermined at junctions, hence we choose to
have some fixed rules for the distribution of traffic plus
an optimization criteria for the flux. We prove existence, uniqueness and
stability of solutions to the Cauchy problem.

Our method is based on wave front tracking approach, see \cite{B},
and works also for boundary data and time dependent coefficients
of traffic distribution at junctions, so including traffic lights.
\end{abstract}
\vspace{0.2cm}

\vspace{2cm}

\begin{center}
Ref. S.I.S.S.A. 13/2002/M
\end{center}

\vfill
\pagebreak

\setlength{\voffset}{-0in}  
\setlength{\textheight}{0.9\textheight}

\thispagestyle{empty}
\null
\vfill
\pagebreak

\thispagestyle{plain}
\setcounter{page}{1}
\setcounter{equation}{0}

\section[]{Introduction} 
\label{section1}
This paper deals with a fluidodynamic model of heavy traffic on a road 
network.
More precisely, we consider the  conservation law formulation proposed by Lighthill and Whitham
\cite{Ligh-W} and Richards \cite{Richards}.
This nonlinear framework is based simply on the conservation of cars and is described by the
equation:
\be
\rho_t+f(\rho)_x=0,\label{eq:conslaw}
\ee
where $\rho=\rho(t,x)\in [0,\rho_{max}]$, $(t,x)\in\R^2$, is the {\em density}
 of cars, $v(t,x)$ is the velocity  and $f(\rho)=v\,\rho$ is the flux. This
model is appropriate to reveal shocks formation as it is natural for
conservation laws, whose solutions may develop discontinuities in finite time
even for smooth initial data, \cite{B}. In most cases one assumes that $v$ is a
function of $\rho$ only and that the corresponding flux is a concave function.
We make the same assumption, moreover we let $f$ have a unique maximum
$\sigma\in ]0,\rho_{max}[$ and for notational simplicity assume $\rho_{max}=1$.

Here we deal with a network of roads, as in \cite{HR}. This means that we
have a finite number of roads modeled by intervals $[a_i,b_i]$ (with one
of the two endpoints possibly infinite) that meet at some junctions. For
endpoints that do not touch a junction (and are not infinite), we assume
to have a given boundary data and solve the corresponding boundary
problem, as in \cite{A, AC, AM, BLN}.  The key role is played by junctions
at which the system is underdetermined even after prescribing the
conservation of cars, that can be written as the Rankine-Hugoniot
relation:
\be
\sum\limits_{i=1}^n f(\rho_i (t, b_i)) = \sum\limits_{j=n
+1}^{n+m}f(\rho_j (t, a_j)),
\label{eq:R-H}
\ee
where $\rho_i$, $i=1,\ldots,n$, are the car densities on incoming roads, while
$\rho_j$, $j=n+1,\ldots,n+m$, are the car densities on outgoing roads.
In \cite{HR}, the Riemann problem, that is the problem with constant initial data on each road,
is solved maximizing a concave function of the fluxes and it is proved existence of weak solutions for
Cauchy problems with suitable initial data of bounded variation.
In this paper we assume that:
\begin{description} 
\item[(A)] there are some prescribed preferences of drivers, that is
the traffic from incoming roads is distributed on outgoing roads according to fixed
coefficients;
\item[(B)] respecting (A), drivers choose so as to maximize fluxes.
\end{description}
To deal with rule (A), we fix a matrix 
$$A\doteq \{\alpha_{ji}\}_{j=n+1,...n+m,\>i=1,...,n}\in \R^{m\times n},$$
such that 
\be 
 \alpha_{ji}\not=\alpha_{ji'}, \qquad 0 < \alpha_{jn} <   1,\qquad
\sum\limits_{j=n +1}^{n+m} \alpha_{ji} =1,\label{matrice}\ee
 for each $i'\not= i=1,...,n$ and $j=n+1,...,n+m,$
where $\alpha_{ji}$ is the percentage of drivers arriving from the
$i-$th incoming road that take the $j-$th outgoing road. Notice that with only
the rule (A) Riemann problems are still underdetermined. This choice
represents a situation in which drivers have a final destination, hence
distribute on outgoing roads according to a fixed law, but maximize the flux
whenever possible. We are able to solve uniquely Riemann problems and, in case
of simple junctions with two incoming and outgoing roads, to generate a
Lipschitz semigroup, defined on $L^1$,  whose trajectories are weak solutions
and respect rules (A) and (B) in case of bounded variation (the same
conditions are not meaningful if the solution is only $L^1$). Our main
technique is the use of a front tracking algorithm and suitable approximations
and functionals to control the total variation. We refer the reader to
\cite{B} for the general theory of conservation laws and for a discussion of wave
front tracking algorithms.

The main difficulty in solving systems of conservation laws is the control
of the total variation, see \cite{B}. It is easy to see that for a single
conservation law the total variation is decreasing, however in our case it
may increase due to interaction of waves with junctions. The problem is
quite delicate, as shown in Appendix B, where an example is given of a
single wave of arbitrarily small strength (variation) that, interacting
with a junction, generates waves whose strengths are bounded away from
zero. Hence we can not expect any bound on the total variation of the
solution in term of the total variation of the initial data, as it is the
case for systems. This arbitrarily large magnification of total variation is
possible only if waves crossing the value $\sigma$ interact with junctions
at which the boundary data of the roads are {\em bad}, that is in
$[0,\sigma]$ for incoming roads and in $[\sigma,1]$ for outgoing roads. We
thus have first to deal with special data of bounded variation that have a
finite number of crossing of the value $\sigma$. The sum of the number of
these crossing plus the number of {\em bad} boundary data is proved to
be decreasing along front tracking approximate solutions.

However, this is not enough since the variation can still increase due to
interactions with junctions (and there is no bound on the number of
interactions).  The conserved quantity is the total variation of the flux. We
prove this fact for junctions with only two incoming roads and two outgoing
ones, and show in Appendix A that the conclusion does not hold for
for junctions with three (or
more) incoming and outgoing roads. 
Unfortunately the total variation of the flux is not
equivalent to the total variation of $\rho$, since $f'(\sigma)=0$. We thus
have to approximate the flux with one having never vanishing derivative and a
corner at $\sigma$, and then pass to the limit.

Our techniques are quite flexible, so we can deal with time dependent
coefficients for the rule (A). In particular we can model traffic lights and
also in this case the control of total variation is extremely delicate. An
arbitrarily small change in the coefficients can produce waves whose strength
is bounded away from zero. Still it is possible to consider periodic
coefficients, a case of particular interest for applications. We can also
deal with roads with different fluxes: this can be treated in the same way
with the necessary notational modifications.

There is an interesting ongoing discussion on hydrodynamic modelization for 
heavy traffic flow. In particular some models using systems of two conservation laws
have been proposed, see \cite{AwR,Colombo,Greenberg}. We do no treat this aspect.

The paper is organized as follows. in Section \ref{section-basic} we give the
definition of weak entrtopic solution and following (A) and (B) we introduce
an admissibility condition. In Section \ref{section2} we prove the existence
and uniqueness of admissible solutions for the Riemann Problem in a junction,
then using this we describe the construction of the approximants for the
Cauchy Problem (see Section \ref{section23}). In section \ref{section3} we
prove the monotonicity of the number of big waves for piecewise constant
solutions. Assuming that $f'$ is bounded away from $0$ and that there are at
most two incoming and outcoming roads in each junction we prove the
monotonicity of the total variation of the flux (see Section \ref{section4})
and existence, uniqueness and stability of admissible solutions for the
Cauchy Problem with suitable $BV$ initial data. Using these results we show
the existence of a unique Lipschitz semigroup of solutions defined on $L^1$
(see Section \ref{section7}) also in the case in which $f$ is smooth. 
In Section \ref{section85} we describe what happens
when there are traffic lights and time dependent coefficients. In
Appendix \ref{section91} we show with an example that the total variation of
the flux can increase when there are three incoming and three outcoming roads
in a junction. Finally, in Appendix \ref{section9} we show that the
interaction of a small wave with a junction can produce a uniformly big
wave. \indent



 
\section[]{Basic Definitions} 
\label{section-basic}

We consider a network of roads, that is a modeled by a finite collection
of intervals $I_i=[a_i,b_i]\subset\R$, $i=1,\ldots,N$, possibly with
either $a_i=-\infty$ or $b_i=+\infty$, on which we consider the equation
(\ref{eq:conslaw}). Hence the datum is given by a finite collection of
functions $\rho_i$ defined on $[0,+\infty[\times I_i$.

On each road $I_i$ we want $\rho_i$ to be a weak entropic solution, 
that is for  $\varphi:I_i\to\R$ smooth with compact support on 
$]0,+\infty[\times ]a_i,b_i[$
\be
\int_0^{+\infty} \int_{a_i}^{b_i} \Big( \rho_i {{\partial \varphi}\over
{\partial t}} + f(\rho_i) {{\partial \varphi}\over
{\partial x}}\Big)dxdt =0, \label{eq:weaksol-oneroad} 
\ee
and for every $k\in\R$ and $\tilde\varphi:I_i\to\R$ smooth, positive with compact support on 
$]0,+\infty[\times ]a_i,b_i[$
\be
\int_0^{+\infty} \int_{a_i}^{b_i} \Big( |\rho_i -k|{{\partial \tilde\varphi}\over
{\partial t}} + \sgn(\rho_i-k)(f(\rho_i)- f(k)){{\partial \tilde\varphi}\over
{\partial x}}\Big)dxdt \geq 0, \label{eq:entsol-oneroad} 
\ee

It is well known that for any initial data in $L^{\infty}$, defined on the
whole $\R$, there exists a unique weak entropic solution depending in a
Lipschitz continuous way from the initial data in $L^1_{loc}$.

We assume that the roads are connected by some junctions.  Each junction
$J$ is given by a finite number of incoming roads and a finite number of
outgoing roads, thus we identify $J$ with 
$((i_1,\ldots,i_n),(j_1,\ldots,j_m))$ where the first $n$--tuple 
indicates the set of incoming roads and the second $m$--tuple indicates 
the set of outgoing roads. 
We assume that each road can be incoming road at most for one
junction and outgoing at most for one junction.

Hence the complete model is given by a couple $(\II,\JJ)$, where
$\II=\{I_i:i=1,\ldots,N\}$ is the collection of roads and
$\JJ$ the collection of junctions.

Fix a junction $J$ with incoming roads, say $I_1,\ldots,I_n$, and outgoing roads, say
$I_{n+1},\ldots,I_{n+m}$. 
A weak solution at the junction $J$ is a collection of
functions $\rho_i:[0,+\infty[\times I_i\to\R$, $i=1,\ldots,n+m$, such that
\be 
\sum\limits_{l=0}^{n+m}\bigg(
\int_0^{+\infty} \int_{a_l}^{b_l} \Big( \rho_l {{\partial \varphi_l}\over
{\partial t}} + f(\rho_l) {{\partial \varphi_l}\over
{\partial x}}\Big)dxdt \Big)=0 \label{weaksol} 
\ee
for each $\varphi_1,...,\varphi_{n+m}$ smooth having compact support in
$]0,+\infty[\times \R$, that are
also smooth across the junction,i.e. 
$$\varphi_i (\cdot,b_i) =\varphi_j (\cdot,a_j), \qquad {{\partial
\varphi_i}\over {\partial x}} (\cdot,b_i)={{\partial
\varphi_j}\over {\partial x}}(\cdot,a_j), \qquad i=1,...,n, \>j=n+1,...,n+m.$$

\begin{remark} 
\label{remark1} Let $ \rho = (\rho_1,\ldots,\rho_{n+m}) $  be a solution of
(\ref{eq:conslaw}) such that each $x\to \rho_i(t,x)$ has bounded variation. We can
deduce that it satisfies the {\it Rankine-Hugoniot Condition} in the junction
$J$, namely 
\be 
\sum\limits_{i=1}^n f(\rho_i (t, b_i)) = \sum\limits_{j=n
+1}^{n+m}f(\rho_j (t, a_j)),
\label{RH}
\ee
for each $t>0.$
\end{remark}

\begin{remark}
The assumption  $\alpha_{ji}\not=\alpha_{ji'}$  in (\ref{matrice})is needed in
order to have uniqueness of solutions to Riemann problems, see Section
\ref{section2}.  However, this condition can be relaxed requiring, for example,
that if $\alpha_{ji}=\alpha_{ji'}$ then the fluxes from road $i$ and $i'$ coincide.
\end{remark}

The rules (A) and (B) can be given explicitly only for solutions with
bounded variation at each time as in next definition.
\begin{definition} \label{definitionsoluz} 

Let $ \rho = (\rho_1,\ldots,\rho_{n+m})$ be such that $\rho_i (t, \cdot)$
is of bounded variation. 
Then $\rho$ is an {\it admissible weak solution} of
(\ref{eq:conslaw}) related to the matrix $A$ satisfying (\ref{matrice}) at the
junction $J$ if and only if following properties hold:
\begin{itemize}
\item[(i)] $\rho$ is a weak solution at the junction;
\item[(ii)]$ f(\rho_j (\cdot, a_j+)) = \sum\limits_{i=1}^n \alpha_{ji}
f(\rho_i (\cdot, b_i-))$, for each $j=n+1,...,n+m$;
\item[(iii)]$ \sum\limits_{i=1}^n f(\rho_i (\cdot, b_i-)) + 
\sum\limits_{j=n+1}^{n+m}f(\rho_j (\cdot, a_j+))$ is maximum subject to (ii).
\end{itemize}
\end{definition}
For every road $I_i=[a_i,b_i]$, if $a_i>-\infty$ and $I_i$ is not the outgoing
of any junction, or $b_i<+\infty$ and $I_i$ is not the incoming road of any
function, then a boundary data $\psi_i:[0,+\infty[\to\R$ is given. 
In this case we ask $\rho_i$ to satisfy
$\rho_i(t,a_i)=\psi_i(t)$ (or $\rho_i(t,b_i)=\psi_i(t)$) in the sense of
\cite{BLN}.

Our aim is to solve the Cauchy problem on $[0,+\infty[$ for a given initial and boundary
data as in next definition. 
\begin{definition}\label{def:sol}
Given $\bar\rho_i:I_i\to\R$ and possibly $\psi_i:[0,+\infty[\to\R$,
functions of $L^{\infty}$, a collection of functions
$\rho=(\rho_1,\ldots,\rho_N)$ with $\rho_i:[0,+\infty[\times I_i\to\R$
continuous as functions from $[0,+\infty[$ into $L^1_{loc}$, is an
admissible solution if $\rho_i$ is a weak entropic solution to
(\ref{eq:conslaw}) on $I_i$, $\rho_i(0,x)=\bar\rho_i(x)$ a.e.,
$\rho_i(t,b_i)=\psi_i(t)$ in the sense of \cite{BLN}, finally such that at
each junction $\rho$ is a weak solution and is an admissible weak solution in
case of bounded variation. \end{definition}
The treatment of boundary data in the sense of \cite{BLN} can be done
in the same way as in \cite{A,AC,AM}, thus we treat the case
without boundary data. All the stated results hold also for
the case with boundary data with the obvious modifications.

On the flux $f$ we make the following assumption
\begin{itemize}
\item[{\bf (${\cal F}$)}] 
$f : [0,1] \rightarrow \R$ is smooth, strictly concave (i.e. $f'' \le -c<0$ for some $c>0$),
$f(0) = f(1) =0$, $\vert f'(x) \vert  \le C < +\infty $ and there exists
$\sigma\in]0,1[$ such that $f'(\sigma)=0$ (that is $\sigma$ is a strict
maximum). \end{itemize}

\section[]{The Riemann Problem} 
\label{section2}

In this section we study Riemann problems. For a scalar conservation law a Riemann problem
is a Cauchy problem for an initial data of Heaviside type, that is piecewise constant
with only one discontinuity. The solutions of these problems are the building
blocks to construct solutions to the Cauchy problem via wave front tracking.
These solutions are formed by continuous waves called rarefactions and traveling discontinuities
called shocks. The speed of waves are related to the values of $f'(\rho)$.

Analogously, we call Riemann problem for the road network the Cauchy
problem corresponding to an initial data that is piecewise constant on each
road. The solutions on each road $I_i$ can be constructed in the same way
as for the scalar conservation law, hence it suffices to describe the
solution at junctions. Because of finite propagation speed, it is
enough to study the Riemann Problem for a single junction.

As explained in the Introduction, we first have to treat the case of fluxes
with nonvanishing derivative, hence we assume that
\begin{itemize}
\item[{\bf (${\cal F}1$)}] 
$f : [0,1] \rightarrow \R$ is continuous, strictly concave, $f(0) =
f(1) =0$ and there exists $\sigma \in]0,1[$ such that $f$ is smooth on 
$[0,\sigma[$ and on $]\sigma, 1]$ and 
\be
0 < c \le \vert f'(x) \vert  \le C <
+\infty, \label{f} \ee
for each $x \in [0, \sigma[ \cup ]\sigma,1].$
\end{itemize}

Consider a junction $J$ in which there are $n$ roads with incoming traffic
and $m$ roads with outgoing traffic. For simplicity we indicate by 
\be
(t,x) \in \, \R_+ \times I_i \mapsto \rho_i (t,x) \in
[0,1],\quad i=1,...,n \label{densentr}
\ee
the densities of the cars on the road with incoming traffic and  
\be
(t,x) \in \, \R_+ \times I_j \mapsto \rho_j (t,x) \in
[0,1],\quad j=n+1,...,n+m \label{densusc}
\ee
those on the roads with outgoing traffic.

\begin{center}
\leavevmode \epsfxsize=2in \epsfbox{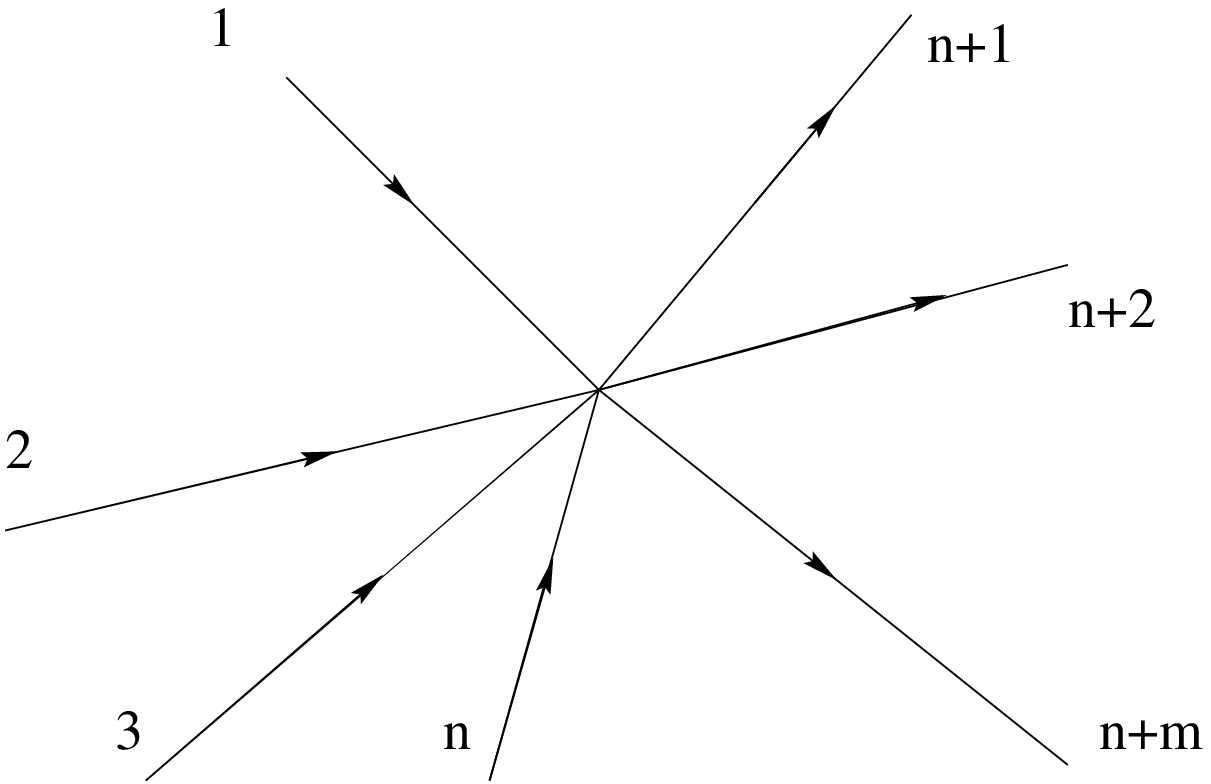}
\end{center}
\begin{center}
figure 1
\end{center} 

We need some more notation:
\begin{definition} \label{deftau}
Let $\tau:[0,1] \rightarrow [0,1]$, $\tau(\sigma)=\sigma$, be the map
satisfying the following
$$\tau(\rho) \not= \rho, \qquad f(\tau(\rho))=f(\rho),$$
for each $\rho \not=\sigma$.
\end{definition}

Clearly $\tau$ is well defined and satisfies
$$0\le \rho \le \sigma \Longleftrightarrow \sigma
\le \tau(\rho) \le 1,\qquad \sigma \le \rho \le 1\Longleftrightarrow 0 \le
\tau(\rho) \le \sigma.$$

The main result of this section is the following Theorem.

\begin{theorem} \label{th1} 
Let $f:[0,1] \rightarrow \R$ satisfy {\bf
(${\cal F}1$)}  and $\rho_{1,0},...,\rho_{n+m,0} \in [0,1]$ be constants. 
There
exists an unique admissible weak solution, in the sense of Definition
\ref{definitionsoluz}, $ \rho =\big(\rho_1,..., \rho_{n+m}\big)$ of
(\ref{eq:conslaw})  at the junction $J$ such that $$ \rho_1 (0, \cdot) \equiv
\rho_{1,0}, ......, \rho_{n+m} (0, \cdot) \equiv \rho_{n+m,0}.$$ Moreover,
there exist a unique $(n+m)-$tuple $(\hat\rho_1,...,\hat\rho_{n+m}) \in
[0,1]^{n+m}$ such that \be \hat \rho_i \in 
\left\{\begin{array}{ll}
\{\rho_{i,0}\} \cup ] \tau(\rho_{i,0}) 1] & \textrm{if $0\le \rho_{i,0}
\le\sigma$},\\  
{[\sigma, 1]} & \textrm{if $\sigma \le \rho_{i,0}\le 1 $},
\end{array} \right. \qquad i=1,...,n \label{ammentr}
\ee 
and 
\be \hat \rho_j \in \left\{\begin{array}{ll}
[0,\sigma] & \textrm{if $0\le\rho_{j,0} \le \sigma$},\\
\{\rho_{j,0}\} \cup [0, \tau(\rho_{j,0})[ & \textrm{if $\sigma \le
\rho_{j,0}\le 1 $}, \end{array} \right. \qquad j=n+1,...,n+m.
\label{ammusc}
\ee 
Fixed $i\in\{1,...,n\}$, if $\rho_{i,0} \le \hat \rho_i$ there results
\be\rho_i(t,x) =
\left\{\begin{array}{ll}
\rho_{i,0} & \textrm{if $x \le {{f(\hat \rho_i)-f(\rho_{i,0})} \over{\hat
\rho_i-\rho_{i,0}}} t$},\\
\hat \rho_i & \textrm{otherwise},
\end{array} \right.
\ee
and if $\hat \rho_i < \rho_{i,0}$
\be\rho_i(t,x) =
\left\{\begin{array}{ll}
\rho_{i,0} & \textrm{if $x \le  f'(\rho_{i,0}) t$},\\
\big(f'\big)^{-1}\big( {x \over t}\big) & \textrm{if $ f'(\rho_{i,0})t\le x \le
f'(\hat\rho_{i}) t$},\\
\hat \rho_i & \textrm{if $ x> f'(\hat\rho_{i})t$}.
\end{array} \right.
\ee
Otherwise, fixed $j\in\{n+1,...,n+m\}$, if $\rho_{j,0} \le \hat \rho_j$ there
results
\be\rho_j(t,x) =
\left\{\begin{array}{ll}
\hat\rho_{j} & \textrm{if $x \le f'(\hat\rho_{j}) t$},\\
\big(f'\big)^{-1}\big( {x \over t}\big) & \textrm{if $ f'(\hat\rho_{j})t\le x
\le f'(\rho_{j,0}) t$},\\
\rho_{j,0} & \textrm{if $ x> f'(\rho_{j,0})t$},
\end{array} \right.
\ee
and if $\hat \rho_j < \rho_{j,0}$
\be\rho_j(t,x) =
\left\{\begin{array}{ll}
\hat\rho_{j} & \textrm{if $x \le {{f(\rho_{j,0})-f(\hat\rho_j)} \over{\hat
\rho_i-\rho_{i,0}}} t$},\\
 \rho_{j,0} & \textrm{otherwise}.
\end{array} \right.
\ee
 \end{theorem}
\vspace{5pt}
\n\textsc{Proof.} Define the map
$$E: (\gamma_1,...,\gamma_n )\in\R^n \longmapsto \sum\limits_{i=1}^n \gamma_i$$
and the sets
$$\Omega_i \doteq 
\left\{\begin{array}{ll}
{[0,f(\rho_{i,0})]} & \textrm{if $0\le \rho_{i,0} \le\sigma $},\\
 {[0, f(\sigma)]} & \textrm{if $\sigma \le \rho_{i,0}\le 1 $},
\end{array} \right. \quad i=1,...,n,$$
$$\Omega_j \doteq 
\left\{\begin{array}{ll}
{[0, f(\sigma)]} & \textrm{if $0\le \rho_{j,0} \le\sigma $},\\
{[0,f(\rho_{j,0})]} & \textrm{if $\sigma \le \rho_{j,0}\le 1 $},
\end{array} \right. \quad j=n+1,...,n+m,$$
$$
\Omega \doteq \Big\{ (\gamma_1,...,\gamma_n )\in \Omega_1 \times .... \times
\Omega_n \big\vert A \cdot (\gamma_1,...,\gamma_n )^T \in \Omega_{n+1} \times
.... \times \Omega_{n+m} \Big\}.
$$ 
Since $E$ is linear, the set $\Omega$ is closed, convex and not empty. 
By (\ref{matrice}) there exists a unique vector
$(\hat\gamma_1,...,\hat\gamma_n )\in \Omega$ such that 
$$E(\hat\gamma_1,...,\hat\gamma_n ) = \max\limits_{(\gamma_1,...,\gamma_n )\in
\Omega} E(\gamma_1,...,\gamma_n ).$$
Fix $i\in\{1,...,n\},$ let $\hat\rho_i \in[0,1]$  be such that
$$f(\hat\rho_i)=\hat\gamma_i, \quad\hat \rho_i \in 
\left\{\begin{array}{ll}
\{\rho_{i,0}\} \cup ] \tau(\rho_{i,0}) 1] & \textrm{if $0\le \rho_{i,0}
\le\sigma$},\\  
{[\sigma, 1]} & \textrm{if $\sigma \le \rho_{i,0}\le 1 $}.
\end{array} \right.
$$
By (${\cal F}1$), $\hat\rho_i $ exists and is unique.
Let
$$\hat\gamma_j \doteq \sum\limits_{i=1}^n \alpha_{ji}\hat\gamma_i, \qquad
j=n+1,...,n+m$$
and $\hat\rho_j \in[0,1]$  be such that 
$$f(\hat\rho_j)=\hat\gamma_j, \quad
\hat \rho_j \in \left\{\begin{array}{ll}
[0,\sigma] & \textrm{if $0\le\rho_{j,0} \le \sigma$},\\
\{\rho_{j,0}\} \cup [0, \tau(\rho_{j,0})[ & \textrm{if $\sigma \le
\rho_{j,0}\le 1 $}. \end{array} \right.$$
Since $(\hat\gamma_1,...,\hat\gamma_n )\in\Omega$,  $\hat\rho_j $ exists and
is unique. Solving the Riemann Problem (see \cite[Chapter 6]{B}) on each road,
the thesis is proved.\hfill$\Box$
\vsp

\section[]{The Wave Front Tracking Algorithm} 
\label{section23}

Once the solution to a Riemann problem is provided, we are able to construct
piecewise constant approximations via wave-front tracking. 
The construction is very similar to that for scalar conservation law, see \cite{B},
hence we only briefly describe it.

Let $\rho_0$ be a piecewise constant map defined on the road network. We want
to construct a solution of (\ref{eq:conslaw}) with initial condition
$\rho(0,\cdot)\equiv \rho_0.$ We begin by solving the Riemann Problems on 
each
road in correspondence of the jumps of $\rho_0$ and the Riemann Problems at
junctions determined by the values of $\rho_0$ (see Theorem \ref{th1}). We
split each rarefaction wave into a rarefaction fan formed by rarefaction
shocks, that are discontinuities traveling with the
Rankine-Hugoniot speed. We always split rarefaction waves inserting the value
$\sigma$ (if it is in the range of the rarefaction), in order to control the
number of big waves defined in next Section.

When a wave interacts with another one we simply solve the new Riemann
Problem, if otherwise it reaches a junction then we solve the Riemann Problem
at the junction. Since the wave speed is bounded there are finitely many
waves on the network at each time $t\ge 0$.  We call the obtained function
{\it a wave front tracking approximate solution}. Given a general initial
data, we approximate it by a sequence of piecewise constant functions and
construct the corresponding approximate solutions. If they converge in
$L^1_{loc}$, then the limit is a weak entropic solution on each road, see
\cite{B} for a proof.

\section[]{Estimates on the Number of Big Waves} 
\label{section3}

In this Section we consider big waves, that are the waves crossing the value
$\sigma$. For these waves the variation of $f(\rho)$ is not comparable to the
variation in $\rho$, more precisely the former can vanish while the second is
different from zero. Since only the variation of $f(\rho)$ happens to be
conserved we need to control the number of big waves.

Let $\rho$ a piecewise constant map defined on the network and $J$ a
junction with $n$ roads with incoming traffic and $m$ roads with
outgoing traffic as in Section \ref{section2}. Define the set
$$
\Phi_J (\rho ) \doteq \big\{ i\in\{1,...,n\}\big\vert
\rho_i(b_i - ) \in [0,\sigma]\big\} \cup \big\{ j\in\{n+1,...,n+m\}\big\vert
\rho_j(a_j + ) \in [\sigma, 1]\big\}.
$$
For each road $I_i$ we denote by 
$\{x_\alpha \in I_i:\alpha\in A_i\}$
and by $\{(\rho(x_\alpha-),\rho(x_\alpha+)):\alpha\in A_i\}$
the set of discontinuity points and the set of discontinuities, respectively,
of the map $\rho_i$ on the road $I_i$. We define
$$
G^i(\rho ) \doteq \big\{ (\rho(x_\alpha-),\rho(x_\alpha+))\big\vert
\ \alpha \in A_i,\ {\rm \sgn} \big(\rho_k(x_\alpha -
)-\sigma\big)\cdot {\rm \sgn}\big(\rho_k(x_\alpha + )-\sigma\big)\leq 0\},$$
with the agreement ${\rm \sgn}(0)=0$, and the functional
$$ N(\rho )\doteq \sum\limits_{J \in\JJ} 
\#\big(\Phi_J (\rho )\big) + \sum_{i=1}^N 
\#\big(G^i(\rho ) \big) ,$$
where $\#$ indicates the cardinality of a set.

The main result of this section is the following.

\begin{lemma} \label{lemma21}
Let $f:[0,1] \rightarrow \R$ satisfy {\bf
(${\cal F}1$)} and $\rho$ be a piecewise constant wave front tracking
approximate solution of (\ref{eq:conslaw}) on the net. Then the map $$t> 0
\longmapsto N\big(\rho(t,\cdot)\big)$$ does not increase.
\end{lemma}

\vspace{5pt}
\n\textsc{Proof.}
We begin considering a single junction $J$ with $n$ roads with incoming
traffic and $m$ roads with outgoing traffic and an equilibrium configuration
$(\rho_{1,0},...,\rho_{n+m,0}) \in [0,1]^{n+m}$, namely the solution of the
Riemann Problem in the junction with that initial data is constant.  Suppose
that a wave on one road arrives to the junction at time $\bar t$ and there is
no other wave on the roads, then we claim that $ N\big(\rho(\bar t
-,\cdot)\big) =N\big(\rho(\bar t +,\cdot)\big)$. Assume that the wave is on
an incoming road, for example the first one and let $(\rho_1, \rho_{1,0})$ be
the values on the left and right side of the wave respectively. Since the
wave is approaching the junction, its speed is positive and so $0 \le \rho_1
\le \sigma$, moreover since it is the unique wave
$$
G^i\big(\rho(\bar t -,\cdot)\big) = \emptyset, \qquad i=2,...,n+m. 
$$  
Let $(\hat\rho_1,...,\hat\rho_{n+m})$ be the solution to the Riemann Problem with
initial data $(\rho_1,\rho_{2,0}...,\rho_{n+m,0})$ (see Theorem \ref{th1}),
there results
$$\hat\rho_1 \in \{\rho_1\}\cup [\sigma,1]$$
$$\hat\rho_i \in \{\rho_{i,0}\}\cup [\sigma,1],\qquad i=2,...,n$$
$$\hat\rho_j \in \{\rho_{j,0}\}\cup [0,\sigma],\qquad j=n+1,...,n+m.$$
In the following we study the change of the functional $N$ due to the presence of new
waves. If a new rarefaction is produced then it can not cross the value $\sigma$, otherwise
there would be rarefaction shocks with positive and negative velocity at the same time.
Hence each functional $G^i$ can not be bigger than one after the interaction.
By abuse of notation, we indicate the whole rarefaction fan as a single 
wave for notational simplicity.

So, fixed $ i\in \{2,...,n\}$, if $i\not\in \Phi_J\big( \rho(\bar t-,
\cdot)\big)$, then
$$i\not\in \Phi_J\big( \rho(\bar t+,\cdot)\big), \qquad
(\rho_{i,0},\hat\rho_i) \not\in G^i\big( \rho(\bar t+,\cdot)\big),$$
and if $i\in \Phi_J\big( \rho(\bar t-,\cdot)\big)$ we have
$$ \hat\rho_i =\rho_{i,0} \Longrightarrow i\in \Phi_J\big( \rho(\bar
t+,\cdot)\big),\>\> (\rho_{i,0},\hat\rho_i) \not\in G^i\big( \rho(\bar
t+,\cdot)\big),$$ 
$$\sigma \le \hat\rho_i \le 1\Longrightarrow i\not\in
\Phi_J\big( \rho(\bar t+,\cdot)\big),\>\> (\rho_{i,0},\hat\rho_i) \in
G^i\big( \rho(\bar t+,\cdot)\big).$$
On the other hand, fixed $ j\in \{n+1,...,n+m\}$, if $j\not\in \Phi_J\big(
\rho(\bar t-, \cdot)\big)$, then
$$j\not\in \Phi_J\big( \rho(\bar t+,\cdot)\big), \qquad
(\hat\rho_j, \rho_{j,0}) \not\in G^j\big( \rho(\bar t+,\cdot)\big),$$
and if $j\in \Phi_J\big( \rho(\bar t-,\cdot)\big)$ we have
$$ \hat\rho_j =\rho_{j,0} \Longrightarrow j\in \Phi_J\big( \rho(\bar
t+,\cdot)\big),\>\> (\hat\rho_j ,\rho_{j,0}) \not\in G^j\big(
\rho(\bar t+,\cdot)\big),$$ 
$$0\le \hat\rho_j \le\sigma \Longrightarrow j\not\in
\Phi_J\big( \rho(\bar t+,\cdot)\big),\>\> (\hat\rho_j, \rho_{j,0})
\in G^j\big( \rho(\bar t+,\cdot)\big).$$
Hence the contribution to $N$ due to roads $I_i$, $i=2,\ldots,n+m$,
does not increase. Let us now treat the waves on the first road.

Notice that if $\rho_{1,0}= \sigma$ then $\rho_1\not=\sigma$ hence
$$1\in \Phi_J\big( \rho(\bar t-,\cdot)\big), \qquad
(\rho_1,\rho_{1,0})\in G^1\big( \rho(\bar t-,\cdot)\big)$$
and $N$ can not increase. The same conclusion holds if
$0\le \rho_{1,0} < \sigma$ and $\rho_1=\sigma$.

Now, if $0\le \rho_{1,0} < \sigma$ and $\rho_1\not=\sigma$, then there results
$$1\in \Phi_J\big( \rho(\bar t-,\cdot)\big), \qquad
(\rho_1,\rho_{1,0})\notin G^1\big( \rho(\bar t-,\cdot)\big)$$
and
$$ \hat\rho_1 =\rho_1 
\Longrightarrow 1\in \Phi_J\big( \rho(\bar
t+,\cdot)\big),\>\> (\rho_1,\hat\rho_1)\not\in G^1\big( \rho(\bar
t+,\cdot)\big),$$
$$ \hat\rho_1 \not=\rho_1\Longrightarrow 
\sigma < \hat\rho_1 \le 1\Longrightarrow 1\not\in
\Phi_J\big( \rho(\bar t+,\cdot)\big),\>\> (\rho_{1},\hat\rho_1) \in
G^1\big( \rho(\bar t+,\cdot)\big).$$
If $\sigma< \rho_{1,0} \le 1$  and $\rho_1 =\sigma $ we have
$$1 \not\in\Phi_J\big( \rho(\bar t-,\cdot)\big), \qquad (\rho_1,
\rho_{1,0})\in G^1\big( \rho(\bar t-,\cdot)\big)  $$
and
$$\hat\rho_1 =\rho_1=\sigma \Longrightarrow 1 \in\Phi_J\big( \rho(\bar
t+,\cdot)\big),\>\> (\rho_1,\hat\rho_1)\not\in G^1\big( \rho(\bar
t+,\cdot)\big),$$
$$\sigma < \hat\rho_1 \le 1\Longrightarrow 1\not\in
\Phi_J\big( \rho(\bar t+,\cdot)\big),\>\> (\rho_{1},\hat\rho_1) \in
G^1\big( \rho(\bar t+,\cdot)\big).$$
Finally, if $\sigma< \rho_{1,0} \le 1$  then $\rho_1 \not=\sigma $, $\hat\rho_1\not=\sigma$, 
and we have
$$
1 \not\in\Phi_J\big( \rho(\bar t-,\cdot)\big), \qquad 
(\rho_{1},\rho_{1,0}) \in G^1\big( \rho(\bar t-,\cdot)\big),
$$ 
hence
$$\hat\rho_1 =\rho_1 \Longrightarrow 1 \in\Phi_J\big( \rho(\bar
t+,\cdot)\big),\>\> (\rho_1,\hat\rho_1)\not\in G^1\big( \rho(\bar
t+,\cdot)\big),$$
$$\sigma < \hat\rho_1 \le 1\Longrightarrow 1\not\in
\Phi_J\big( \rho(\bar t+,\cdot)\big),\>\> (\rho_{1},\hat\rho_1) \in
G^1\big( \rho(\bar t+,\cdot)\big).$$
We conclude
$$N\big( \rho(\bar t+,\cdot)\big) = N\big( \rho(\bar t-,\cdot)\big).$$ 

The conclusion can be obtained in the same way if the wave is arriving to the
junction from another road.
Moreover, when two waves interact on the same road then there is a cancellation
or gluing of waves and it is easy to check that $N$ is constant or decreases.
The proof is concluded.
\hfill$\Box$
\vsp

\section[]{Estimates on Flux Variation} 
\label{section4}
\indent

This Section is dedicated to the estimation of the total variation of the flux
along a solution. We assume that every junction has at most two incoming 
roads
and two outgoing ones. This hypothesis is crucial, because, as shown in the
Appendix \ref{section91}, the presence of more complicate junctions provokes
increase of the total variation of the flux.

\begin{lemma} 
\label{lemma31} Let $f:[0,1] \rightarrow \R$ satisfy {\bf
(${\cal F}1$)}. Consider a network $(\II,\JJ)$ in which each junction has at 
most
two incoming road and two outgoing ones. Let $\rho$ be a piecewise constant
wave front tracking  approximate solution.Then the map  
$$
t> 0\,\mapsto\, \tv \big( f(\rho( t, \cdot))\big) , 
$$ 
is not increasing.
\end{lemma}

\vspace{5pt}
\n\textsc{Proof.} First of all we consider a single junction $J$ with $n\le2$
roads with incoming traffic and $m\le2$ roads with outgoing traffic as in
Section \ref{section2}. It suffices to study the case $n=m=2$, the other ones
are simpler. Let $(\rho_{1,0},...,\rho_{4,0})$ be  an equilibrium
configuration in the junction $J$. Assume that a wave comes to the junction 
at
the time $\bar t$, we claim that  \be \tv \big( f(\rho( \bar t +, \cdot))\big)
= \tv \big( f(\rho( \bar t -, \cdot))\big). \label{flux} \ee
We begin assuming that the wave is on an incoming road, for example the first
one, and that it is given by the values $(\rho_1, \rho_{1,0}).$ 
Let us define the incoming flux
 \be
f^{in}(y) \doteq 
\left\{\begin{array}{ll}
{f(y)} & \textrm{if $0\le y \le\sigma $},\\
 {f(\sigma)} & \textrm{if $\sigma \le y\le 1 $},
\end{array} \right. \label{influx}\ee
and the outgoing flux
\be
f^{out}(y) \doteq 
\left\{\begin{array}{ll}
{f(\sigma)} & \textrm{if $0\le y \le\sigma $},\\
 {f(y)} & \textrm{if $\sigma \le y\le 1 $}.
\end{array} \right. \label{outflux}\ee
Clearly, since the wave on the first road has positive velocity, we have
\be 
0\le \rho_1 \le \sigma, \qquad f(\rho_1)< f^{out}(\rho_{1,0}).\label{ro1}\ee
Let $(\hat\rho_1,...,\hat\rho_{4})$ the solution of the Riemann Problem
in the junction $J$ with initial data $(\rho_1,...,\rho_{4,0})$
(see Theorem \ref{th1}). By definition $\big(
f(\rho_{1,0}),f(\rho_{2,0})\big)$ is the maximum of the map $E$ on the
domain 
$$\Omega_0 \doteq \Big\{ (\gamma_1,\gamma_2 )\in \Omega_{1,0}
\times  \Omega_{2,0}\big\vert A \cdot (\gamma_1,\gamma_2 )^T \in
 \Omega_{3,0}\times \Omega_{4,0}\Big\},$$ 
and $\big( f(\hat \rho_{1}),f(\hat \rho_{2})\big)$ is the maximum of
the map $E$ on the domain
$$\hat\Omega \doteq \Big\{ (\gamma_1,\gamma_2 )\in \Omega_1
\times\Omega_{2,0} \big\vert A \cdot
(\gamma_1,\gamma_2 )^T \in \Omega_{3,0}\times \Omega_{4,0}\Big\},$$
where
$$\Omega_{j,0} \doteq \left\{\begin{array}{ll}
{[0, f^{in}(\rho_{j,0})]} & \textrm{if $j=1,2 $},\\
 {[0, f^{out}(\rho_{j,0})]} & \textrm{if $j=3,4$},
\end{array} \right.$$
and ,by (\ref{ro1}),
$$\Omega_{1} \doteq [0, f^{in}(\rho_1)]= [0, f(\rho_1)].$$
It is also clear that 
$$\big( f(\rho_{1,0}),f(\rho_{2,0})\big) \in \partial \Omega_0, \qquad
\big( f(\hat \rho_{1,0}),f(\hat \rho_{2,0})\big)\in \partial \hat\Omega.$$
To simplify the notations, define
$$\alpha_1\doteq\alpha_{31},\qquad \alpha_2\doteq\alpha_{32}$$
then, by (\ref{matrice}),
$$1-\alpha_1 = \alpha_{41},\qquad 1-\alpha_2 = \alpha_{42}$$
We distinguish two cases. First we suppose that 
\be
f(\rho_1)< f(\rho_{1,0}),\label{case1}
\ee
(equality can not happen in the previous equation because the wave would have velocity zero).
Then there results $\hat\Omega\subset \,\Omega_0$, hence
\be f(\hat\rho_{1}) \le f(\rho_{1})< f(\rho_{1,0}),\qquad
f(\hat\rho_{1})+f(\hat\rho_2) \le f(\rho_{1,0})+f(\rho_{2,0}),\label{fl1}\ee
where the first inequality is due to the fact that the wave $(\rho_1,\hat\rho_1)$ has
negative velocity. We claim that
\be f(\rho_{2,0}) \le f(\hat\rho_2), \label{fl2}\ee
and
\be  f(\hat\rho_3)\le f(\rho_{3,0}), \qquad 
f(\hat\rho_4)\le f(\rho_{4,0}).\label{fl3}\ee  
The points $\big( f(\rho_{1,0}),f(\rho_{2,0})\big),\>
\big( f(\hat\rho_{1,0}),f(\hat \rho_{2,0})\big)$ 
are on the boundaries of $\Omega_0$, $\hat\Omega$ respectively,
where $E$ is maximum, hence they are on one of the curves
$$\alpha_{1}\gamma_1+ \alpha_{2}\gamma_2= f^{out}(\rho_{3,0}),\quad
\quad  (1-\alpha_{1})\gamma_1+ (1-\alpha_{2})\gamma_2=
f^{out}(\rho_{4,0}),\quad \quad \gamma_2 = f^{in}(\rho_{2,0}).$$ 
Using (\ref{case1}), we immediately get (\ref{fl2}).  Let us assume that the two points
are on the same curve, the general case being similar, for example on
$$
\alpha_{1}\gamma_1+\alpha_{2}\gamma_2= f^{out}(\rho_{3,0}).
$$
From (\ref{case1}) it follows that the map $E$ is increasing on the curve 
$$
\gamma_1 \mapsto\Big(\gamma_1,\, f^{out}(\rho_{3,0})-{{\alpha_{1}}\over{\alpha_{2}}}\gamma_1\Big),
$$ 
otherwise we contradict the maximality of $E$ at $\big( f(\rho_{1,0}),f(\rho_{2,0})\big)$. Thus
$\alpha_1< \alpha_2$, $\hat\rho_1=\rho_1$, and 
$$ 
f(\hat\rho_{1})= f(\rho_{1}),\qquad
f(\hat\rho_{3})=f(\rho_{3,0})=f^{out}(\rho_{3,0}).
$$ 
On the other hand, by (\ref{fl1}) and (\ref{fl2}), we have 
$$
f(\hat\rho_{4})=(1-\alpha_{1})f(\hat\rho_1)+(1-\alpha_{2})f(\hat\rho_2)\le
$$
$$
\le(1-\alpha_{1})\big(f(\rho_{1,0})+f(\rho_{2,0})-f(\hat\rho_2)\big)+
(1-\alpha_{2})f(\hat\rho_2)=
$$
$$
=(1-\alpha_{1})\big(f(\rho_{1,0})+f(\rho_{2,0})\big)+\big(\alpha_1-
\alpha_2\big)f(\hat\rho_2)\le
$$
$$\le
(1-\alpha_{1})\big(f(\rho_{1,0})+f(\rho_{2,0})\big)+\big(\alpha_1-
\alpha_2\big)f(\rho_{2,0})=f(\rho_{4,0}).
$$ 
Using the Rankine Hogoniot Condition (\ref{RH}) in
the junction  (\ref{fl1}), (\ref{fl2}) and (\ref{fl3}), we get   
$$
\tv\big( f(\rho( \bar t +, \cdot))\big) =   \vert f(\hat\rho_1) -
f(\rho_1)\vert + \vert f(\hat\rho_2) - f(\rho_{2,0})\vert +\vert f(\hat\rho_3)
- f(\rho_{3,0})\vert + \vert f(\hat\rho_4) - f(\rho_{4,0})\vert=
$$ 
$$
=\big(f(\hat\rho_2)-f(\rho_{2,0}) \big)+  \big(f(\rho_{3,0})-f(\hat\rho_3)\big)+ \big(f(\rho_{4,0}) -
f(\hat\rho_4)\big)=
$$
$$
=f(\rho_{1,0})- f(\hat\rho_1)   =f(\rho_{1,0})- f(\rho_{1})
= \tv \big( f(\rho( \bar t -,\cdot))\big). 
$$
Suppose now that 
$$f(\rho_{1,0})< f(\rho_1),$$
then $\Omega_0\subset \hat\Omega$
and using the previous arguments
$$
f(\hat\rho_1) = f(\rho_{1}), \quad f(\hat\rho_2)\le
f(\rho_{2,0}), \quad f(\rho_{3,0}) \le  f(\hat\rho_3), \quad f(\rho_{4,0}) \le
f(\hat\rho_4).
$$ 
By the Rankine Hogonot Condition in the junction (see (\ref{RH})),  we have 
$$
\tv\big( f(\rho( \bar t +, \cdot))\big) = 
\vert f(\hat\rho_1) - f(\rho_1)\vert +
\vert f(\hat\rho_2) - f(\rho_{2,0})\vert +\vert f(\hat\rho_3) -
f(\rho_{3,0})\vert + \vert f(\hat\rho_4) - f(\rho_{4,0})\vert=
$$ 
$$
=\big(f(\rho_{2,0})-f(\hat\rho_{2}) \big)+ 
\big(f(\hat\rho_{3})-f(\rho_{3,0})\big)+
\big(f(\hat\rho_{4})-f(\rho_{4,0}))\big)=
$$ 
$$
= \big(f(\hat\rho_1)-f(\rho_{1,0}))\big) + 
=f(\rho_{1})-f(\rho_{1,0})= \tv \big( f(\rho( \bar t -,\cdot))\big). 
$$
In the general case we have only to observe that the total variation of $\rho$
does not increase on the roads (see \cite[Chapter 6]{B}) and when a wave
approaches a junction we can use the previous argument, so  the proof is
concluded.\hfill$\Box$
\vsp

\section[]{Solutions with a finite number of Big Waves} 
\label{section5}
\indent

In this section we prove existence and stability of solutions with a finite
number of big waves.

\begin{definition}
We call ${\cal D}_n $ the set of all maps $\bar \rho: \cup_i I_i\mapsto \R$ 
defined on the network with bounded variation such that there exists a 
sequence 
$\{\bar\rho_\nu\}$ of piecewise constant maps satisfying
\be
\tv (\bar\rho_\nu) \le \tv(\bar\rho), \quad\quad \> N(\bar\rho_\nu) \le n, 
 \label{dn1}\ee
for each $\nu\in\NN $ and 
\be
\bar\rho_\nu \longrightarrow \bar\rho \quad{\rm in }\> L^1.
\label{dn2}\ee
\end{definition}
Notice that if $\bar\rho\in{\cal D}_n$ then
\be\tv(\bar\rho) \le \Big \Vert {1\over {f'}}\Big \Vert_{L^{\infty}} \tv\big( f(\bar\rho) \big)
+n. \label{tv}
\ee
The existence of solutions with values in the domain ${\cal D}_n$ is ensured
by next Theorem.
\begin{theorem} 
\label{existence}
 Let $f:[0,1] \rightarrow \R$ satisfying {\bf (${\cal F}1$)}. Consider a road
network $(\II,\JJ)$ in which all junctions have at most two incoming roads
and two outgoing ones. Given $n \in \NN$ and $\bar \rho \in {\cal D}_n$,
there exists an admissible solution $\rho$ in the sense of Definition
\ref{def:sol} such that $\rho (0, \cdot) = \bar\rho$, and $\rho(t, \cdot) \in
{\cal D}_n,$ for each $t\ge 0$. \end{theorem}
\vspace{5pt}
\n\textsc{Proof.} 
Let $\{\bar\rho_\nu\} \subset {\cal D}_n$ be a sequence of
piecewise constant maps such that 
$$
\tv (\bar\rho_\nu) \le \tv (\bar\rho),
\qquad \bar\rho_\nu \longrightarrow \bar \rho \> \>{\rm in} \>\> L^1,
$$ 
and $\rho_\nu$ the wave front tracking approximate solutions such that 
$\rho_\nu (0, \cdot) = \bar\rho_\nu.$ By Lemma \ref{lemma21} we have 
$$ 
\rho_\nu (t,\cdot) \in {\cal D}_n, \qquad t\ge 0, \>\> \nu\in \NN 
$$ 
and, by (\ref{tv}) and Lemma \ref{lemma31}, 
$$ 
\tv (\rho_\nu (t, \cdot)) \le \Big \Vert {1\over
{f'}}\Big \Vert_{L^{\infty}} \tv\big( f(\rho_\nu (t, \cdot)) \big) +n\le 
$$
$$ 
\le \Big \Vert {1\over {f'}}\Big \Vert_{L^{\infty}} \tv\big( f( \bar
\rho_\nu) \big)+n \le \Big \Vert {1\over {f'}}\Big \Vert_{L^{\infty}} \cdot
\Vert f'\Vert_{L^{\infty}}\tv(\bar \rho) +n.  
$$ 
Since the wave speeds are bounded, the maps $\rho_\nu$ are uniformly
Lipschitz continuous from $[0,+\infty[$ into $L^1_{loc}$ for the $L^1$ norm
on every compact set, and are obviously uniformly bounded. Then by Helly
Theorem (see \cite[Theorem 2.4]{B}), $\rho_\nu$ converge to some
continuous map $\rho\in L^1_{loc}([0,+\infty[\times \cup_i I_i,\R)$ such
that, up to redefining $\rho$ on a set of zero measure, $\rho(t,\cdot)$ has
bounded variation. Since for every $t$ we can obtain $\rho(t,\cdot)$ as
limits of $\rho_\nu(t_n,\cdot)$ for some $t_n\to t$, we get
$\rho(t,\cdot)\in{\cal D}_n$ for every $t\geq 0$.

It is a standard argument, see \cite{B}, to prove that $\rho$ solves the
conservation law on the road network. Moreover, $\rho(t,\cdot)\in BV$ and one
can easily check the other properties to guarantee that $\rho$ is an
admissible solution. The proof is concluded.\hfill$\Box$
\vspace{5pt}

\n Regarding stability of solutions we have:
\begin{theorem} \label{stability} Let $f:[0,1] \rightarrow \R$ satisfy {\bf
(${\cal F}1$)}. Consider a network in which all junctions have at most two
incoming roads and two outgoing ones.  Let $n \in \NN,$
$\rho$ and $ \tilde \rho$ be admissible solutions such that
$$\rho (t, \cdot), \tilde\rho (t, \cdot) \in {\cal D}_n\cap L^1, \qquad t\ge 0.$$
There results
$$\Vert \rho (t, \cdot)- \tilde\rho (t, \cdot)\Vert_{L^1} \le 
\Vert \rho (0, \cdot)- \tilde\rho (0, \cdot)\Vert_{L^1},$$
for each $t\ge 0$.
 \end{theorem}
We begin proving a lemma.
\begin{lemma} 
\label{lemma41} Let $f:[0,1] \rightarrow \R$ satisfy {\bf
(${\cal F}1$)} and $J$ be a junction with two incoming roads and two outgoing
ones as in Section \ref{section2}. 
Let $\rho_1, \rho_{1,0},..., \rho_{4,0}
\in [0,1], \> \xi \in\R$ and  $\bar x < b_1$ be such that
$$\bar x +\xi < b_1  , \qquad 
\rho_1 \in \left\{\begin{array}{ll}
{[0, \sigma]\setminus \{\rho_{1,0}\}} & \textrm{if $0\le \rho_{1,0} \le
\sigma$},\\
 {[0, \tau(\rho_{1,0})]} & \textrm{if $\sigma \le \rho_{1,0} \le 1$}.
\end{array} \right.
$$
Let $\rho$ and $\rho^{\xi}$ be the wave front tracking approximate solutions 
of (\ref{eq:conslaw}) on
$J$ such that
$$\rho_1(0, \cdot) = \rho_1 \cdot \chi_{]-\infty, \bar x]} + \rho_{1,0} \cdot
\chi_{]\bar x, b_1]}, \quad \rho_1^{\xi}(0, \cdot) = \rho_1 \cdot
\chi_{]-\infty, \bar x+\xi]} + \rho_{1,0} \cdot \chi_{]\bar x+\xi, b_1]},$$
$$\rho_i (0, \cdot) \equiv \rho_i^{\xi}(0, \cdot)\equiv\rho_{i,0}, \qquad
i=2,3,4.$$
There results
$$\Vert \rho_1^{\xi}(t, \cdot) -\rho_1(t, \cdot)\Vert_{L^1}\le \vert \xi \vert
\cdot \vert \rho_{1,0} - \rho_{1}\vert = \Vert \rho_1^{\xi}(0, \cdot)
-\rho_1(0, \cdot)\Vert_{L^1},$$
for each $t\ge 0.$
\end{lemma}

\vspace{5pt}
\n\textsc{Proof.} 
We begin assuming that $(\rho_{1,0},..., \rho_{4,0})$ is an equilibrium
configuration. By possibly  changing the notations, we can assume that $\xi
>0.$ Since we approximate the rarefaction fronts with many small shocks we have
only to study the case in which  the Riemann Problem $(\rho_1, \rho_{1,0})$ on
the first road and the one $(\rho_{1,0},..., \rho_{4,0})$ in the junction 
generate
only shocks. Let $(\hat\rho_1,...,\hat\rho_{4})$ be the solution of the Riemann
Problem with initial data $(\rho_1,\rho_{2,0}...,\rho_{4,0})$ (see Theorem
\ref{th1}) and $\lambda_{1,0}, \hat\lambda_1, \hat\lambda_2, \hat\lambda_3,
\hat\lambda_{4}$ the velocities of the shocks generated by the Riemann
Problems $(\rho_1, \rho_{1,0}), (\rho_1, \hat\rho_{1}),(\rho_2,\hat\rho_{2}),
(\hat \rho_{3}, \rho_{3,0}), (\hat \rho_{4},\rho_{4,0})$ respectively.\\

If $t\ge \displaystyle {{b_1 - \bar x -\xi} \over {\lambda_{1,0}}}$ by the
first part of the proof of Lemma \ref{lemma31} and the Rankine Hugoniot
Conditions on the roads (see \cite[Lemma 4.2]{B}) we have 
$$\Vert
\rho_1^{\xi}(t, \cdot) -\rho_1(t, \cdot)\Vert_{L^1}= {{\vert \xi \vert}\over
{\vert \lambda_{1,0} \vert}} \Big( \vert\rho_{1} - \hat \rho_1 \vert \cdot
\vert \hat\lambda_1\vert + \sum\limits_{i=2}^{4} \vert \rho_{i,0} - \hat
\rho_i \vert \cdot \vert \hat\lambda_i\vert \Big)=$$
 $$= {{\vert \xi
\vert}\over {\vert \lambda_{1,0} \vert}} \Big(\vert f(\rho_{1}) - f(\hat
\rho_1) \vert  + \sum\limits_{i=2}^{4} \vert f(\rho_{i,0})
- f(\hat \rho_i)\vert  \Big)={{\vert \xi
\vert}\over {\vert \lambda_{1,0} \vert}} \tv \big(f(\rho(t,\cdot))\big) = $$
$$={{\vert \xi \vert}\over {\vert \lambda_{1,0} \vert}} \tv
\big(f(\rho(0,\cdot))\big) = {{\vert \xi \vert}\over {\vert \lambda_{1,0}
\vert}}\big\vert f(\rho_1) - f(\rho_{1,0})\big\vert = \vert \xi \vert
\cdot \vert \rho_1 - \rho_{1,0}\vert = \Vert \rho_1^{\xi}(0, \cdot) -\rho_1(0,
\cdot)\Vert_{L^1}.$$

In the case in which $(\rho_{1,0},..., \rho_{4,0})$ is not an equilibrium
configuration we have only to recall that the $L^1-$distance between the
solutions decreases on each roads (see of \cite[Corollary 6.1]{B}) and use the
same arguments.

This concludes the proof.\hfill$\Box$
\vsp

\n\textsc{Proof of Theorem \ref{stability}.} 
Let $\rho^k$ and $\tilde\rho^k$ be front tracking approximate solutions such that 
$$\Vert\rho^k (0,\cdot)-\rho(0,\cdot)\Vert_{L^1} \le {1 \over k},
\qquad \Vert\tilde\rho^k (0,\cdot)-\tilde\rho(0,\cdot)\Vert_{L^1} \le {1 \over k},
\qquad k\in\NN,$$
$$\tv (\rho^k (0,\cdot)) \le \tv (\rho (0,\cdot)),\quad \tv (\tilde\rho^k
(0,\cdot)) \le \tv (\tilde\rho (0,\cdot)), \quad k\in\NN.$$
Now consider finitely many wave front tracking approximate solutions 
$\rho^{k,0},..., \rho^{k,N}$
$$\rho^{k,0}\equiv \rho^{k},\qquad \rho^{k,N}\equiv \tilde\rho^{k}$$
where $\rho^{k,h}$ is obtained by $\rho^{k,h-1}$ shifting and rescaling only 
one
jump as in \cite{BCP} and \cite{BM}. Precisely denoting
$$\rho^{k,h-1}_i (0, \cdot)= \sum\limits_{l=0}^{N-1} \alpha_l \cdot
\chi_{[\beta_l, \beta_{l+1}[},$$
there exist $\lambda, \xi \in\R$ and $\bar l \in \{1,...,N-1\}$ such that
$$\rho^{k,h}_i (0, \cdot)= \sum\limits_{l=0}^{\bar l -2} \alpha_l \cdot
\chi_{[\beta_l, \beta_{l+1}[} 
+\alpha_{\bar l -1} \cdot \chi_{[\beta_{\bar l -1},\beta_{\bar l} + \xi[} +
\lambda \cdot \chi_{[\beta_{\bar l} + \xi, \beta_{\bar
l+1}[} +\sum\limits_{l=\bar l +1}^{N+1} \alpha_l \cdot \chi_{[\beta_l,
\beta_{l+1}[} ,$$
with
$$\beta_{\bar l} -\beta_{\bar l -1} \le \xi \le  \beta_{\bar l+1} -\beta_{\bar
l },$$
for each $i=1,...,N$.
In this way we have 
\be\sum\limits_{h=1}^{N}
\Vert\rho^{k,h}(0,\cdot)-\rho^{k,h+1}(0,\cdot)\Vert_{L^1} = 
\Vert\rho^{k}(0,\cdot)-\tilde\rho^{k}(0,\cdot)\Vert_{L^1}.\label{*}
\ee
Since the distance between solutions decreases on
each road (see \cite[Corollary 6.1]{B}) and  by the previous lemma, we have
$$\Vert\rho^{k,h}(t,\cdot)-\rho^{k,h+1}(t,\cdot)\Vert_{L^1} \le
\Vert\rho^{k,h}(0,\cdot)-\rho^{k,h+1}(0,\cdot)\Vert_{L^1}.$$ So, by (\ref{*}),
$$\Vert\rho^{k}(t,\cdot)-\tilde\rho^{k}(t,\cdot)\Vert_{L^1} \le
\sum\limits_{h=1}^{N}
\Vert\rho^{k,h}(t,\cdot)-\rho^{k,h+1}(t,\cdot)\Vert_{L^1} \le$$
$$
\le\sum\limits_{h=1}^{N}
\Vert\rho^{k,h}(0,\cdot)-\rho^{k,h+1}(0,\cdot)\Vert_{L^1} = 
\Vert\rho^{k}(0,\cdot)-\tilde\rho^{k}(0,\cdot)\Vert_{L^1}.
$$ 
Moreover there exists a decreasing sequence $\{k_n\}\subset\,\NN$ such that 
$\rho^{k_n}\rightarrow \rho$ and $\tilde\rho^{k_n}\rightarrow \tilde\rho$
in $ L^1$ as $k_n\rightarrow +\infty$. Hence
$$\Vert\rho(t,\cdot)-\tilde\rho(t,\cdot)\Vert_{L^1}
\le\Vert\rho(0,\cdot)-\tilde\rho(0,\cdot)\Vert_{L^1},$$
as to be proved.\hfill$\Box$

\section[]{Existence and Stability of Solutions in $L^1$} 
\label{section7}

\indent
Let us first consider the case in which {\bf (${\cal F}1$)} holds true.
\begin{theorem} \label{exstl1} Let $f:[0,1] \rightarrow \R$ satisfy {\bf
(${\cal F}1$)}. Consider a road network in which each junction has at most
two incoming roads and two outgoing ones. Let $\bar \rho$ be an initial data
in $L^1_{loc}$ and fix $T>0$. Then there exists a unique admissible solution
$\rho$ defined on $[0,T]$, obtained as limit of wave front tracking
approximate solutions such that $\rho(0,\cdot)=\bar\rho$. Moreover if $\bar
\rho\in L^1$ then $\rho(t,\cdot)\in L^1$.

If $\rho$ and $\tilde \rho$ are admissible solutions
obtained as limit of wave front tracking approximate solutions such that $\rho(t, \cdot)
,\,\tilde\rho(t, \cdot)\in L^1$, for every  $t\ge 0,$ then for each $t\ge0$
\be
\Vert\rho(t, \cdot) -\tilde\rho(t, \cdot)\Vert_{L^1} \le  
\Vert\rho(0, \cdot) -\tilde\rho(0, \cdot)\Vert_{L^1}.\label{st}
\ee
\end{theorem}
\vsp

\n\textsc{Proof.} We begin proving the existence of a solution for
$\bar\rho\in L^1$.  There exists $\{\bar\rho_n\}$ sequence of piecewise
constant maps defined on the network such that
\be
\bar \rho_n \in {\cal D}_n\cap L^1, \qquad \bar\rho_n \longrightarrow \bar\rho \quad
{\rm in} \>\> L^1. \label{appr}\ee
Let $\rho_n$ be the wave front tracking approximate solutions with $\rho_n 
(0,\cdot) = \bar \rho_n.$
Fix $n < m$, by Lemma \ref{lemma21}, there results
$\rho_n (t,\cdot) \in {\cal D}_n \subset \,{\cal D}_m, \rho_m (t,\cdot)\in
{\cal D}_m$ and by Theorem \ref{stability},
$$
\Vert\rho_n (t,\cdot)- \rho_m (t,\cdot)\Vert_{L^1} \le \Vert\bar\rho_n
- \bar\rho_m \Vert_{L^1}. 
$$
Hence $\{\rho_n (t,\cdot)\}$ is a Cauchy sequence in $L^1\big([0,T]\times\R,\R\big)$. 
Then there exists $\rho$ such that $\rho_n (t,\cdot)\longrightarrow \rho (t,\cdot)$ in 
$L^1\big([0,T]\times\R,\R\big)$. It is easy to check that $\rho$ is as admissible solution.
The case of $L^1_{loc}$ can be obtained by localization.

Now we prove (\ref{st}). Let $\{\rho_n\},\>\{\tilde\rho_n\}$ be sequences of
wave front tracking approximate solutions such that
$\rho_n(t, \cdot),\>\tilde\rho_n(t, \cdot)\in {\cal D}_n\cap L^1$ and
$\rho_n\longrightarrow \rho$, $\tilde\rho_n\longrightarrow\tilde\rho$ in
$L^1\big([0,T]\times\R,\R\big)$. By Theorem \ref{stability}, we have 
$$\Vert\rho_n(t, \cdot) -\tilde\rho_n(t, \cdot)\Vert_{L^1} \le  
\Vert\rho_n(0, \cdot) -\tilde\rho_n(0, \cdot)\Vert_{L^1}.$$
Therefore (\ref{st}) is proved and uniqueness holds true. 
\hfill$\Box$

We now relax the assumption (${\cal F}1$), namely we
suppose that  $f$ satisfies (${\cal F}$). 

Let $\{f_\nu\}$  be a sequence of maps satisfying (${\cal F}1$) such
that 
\be f_\nu (\sigma) = \max\limits_{\rho \in [0,1]} f_\nu (\rho), \qquad \nu\in \NN
 \label{appr2}\ee
and 
\be f_\nu\longrightarrow f\>\> {\rm in} \> L^{\infty}([0,1])\quad{\rm
and}\quad {f_\nu}'\longrightarrow {f}'\>\> {\rm in} \>
L^{\infty}([0,1]). \label{appr3}\ee
Moreover let $\bar\rho$ be an initial data in $L^1_{loc}$. We know that 
there exists a unique $\rho_\nu =\rho_\nu (t,x)$ admissible solution to the 
Cauchy Problem on the network  (see Theorem \ref{exstl1}) obtained as limit
of front tracking approximate solutions for
\be
\rho_t + f_\nu (\rho)_x =0, \qquad\qquad \rho(0,\cdot) \equiv
\bar\rho, \qquad\qquad \nu \in \NN.\label{Cauchy}
\ee
\begin{theorem} \label{exstsmooth}
Let $f:[0,1] \rightarrow \R$ satisfy {\bf (${\cal F}$)}. Consider a road
network in which all the junction have at most two incoming roads and two
outgoing ones. Let $\bar \rho$ be an initial data in $L^1_{loc}$ and fix
$T>0$. Then there exists a unique admissible solution $\rho$ defined on 
$[0,T]$, with $\rho(0,\cdot)=\bar\rho$,
obtained as limit of solutions to (\ref{Cauchy}).
The limit does not depend on the choice of the functions $f_\nu$ and
if $\bar \rho\in L^1$ then $\rho(t,\cdot)\in L^1$.

If $\rho$ and $\tilde \rho$ are such admissible solutions
and satisfy $\rho(t, \cdot),\,\tilde\rho(t, \cdot)\in L^1$, for every  $t\ge 0,$ 
then
\be
\Vert\rho(t, \cdot) -\tilde\rho(t, \cdot)\Vert_{L^1} \le  
\Vert\rho(0, \cdot) -\tilde\rho(0, \cdot)\Vert_{L^1}.\label{sts}
\ee
\end{theorem}
\vsp
The Theorem can be proved exactly as Theorem \ref{exstl1} from next Lemmas.

\begin{lemma} 
\label{lemma51}Let $f:[0,1] \rightarrow \R$ satisfy {\bf
(${\cal F}$)}. Consider a road network in which all junctions have at most
two incoming roads and two outgoing ones. Let $\bar \rho$ be an initial data
in ${\cal D}_n \cap L^1$ and fix $T>0$. Then there exists a unique admissible
solution $\rho$ defined on $[0,T]$, with $\rho(0,\cdot)=\bar\rho$,
obtained as limit of solutions to (\ref{Cauchy}).
The limit does not depend on the choice of the functions $f_\nu$ and
and $\rho(t,\cdot)\in L^1$. 
\end{lemma}

\vspace{5pt}
\n\textsc{Proof.} Let $\{f_\nu\}$ be a sequence of maps satisfying (${\cal
F}1$), (\ref{appr2}) and (\ref{appr3}) and $\rho_\nu$  be the  admissible
solutions for the Cauchy problems associated to $f_\nu$. By Theorem
\ref{existence} we have
\be \rho_\nu(t,\cdot) \in {\cal D}_n, \qquad\qquad \nu\in\NN,\quad 0\le t\le
T.\label{stefano}\ee
Moreover there results
\be\label{pippo}
\Vert\rho_\nu (t, \cdot) -\rho_\mu(t, \cdot)\Vert_{L^1} \le C \cdot
\Vert {f_\nu}' -{f_\mu}'\Vert_{L^{\infty}}\cdot \tv ( {\bar\rho}),
\ee
where $C$ depends only on $f$. 
By (\ref{appr3}), $\{\rho_\nu\}$ is a
Cauchy sequence in $L^1$, then there exists $\rho\in L^1$
such that $\rho_\nu\rightarrow \rho$ in $L^1$. Moreover, $\rho$
is an admissible solution and satisfies $\rho(0, \cdot)\equiv\bar\rho$.
From (\ref{pippo})  we have that $\rho$ does not depend on
the choice of $\{f_\nu\}$, so we are done.\hfill$\Box$

\vsp
\begin{lemma} 
\label{lemma61} Let $f:[0,1] \rightarrow \R$ satisfy {\bf
(${\cal F}$)}. Consider a road network in which all junctions have at most
two incoming roads and two outgoing ones. Fix $T>0$ and let $\tilde \rho,\>
\rho$ be  admissible solutions  in $L^1$, 
obtained as limit of solutions to (\ref{Cauchy}), defined on $[0,T]$. 
If $\rho(0,\cdot),\> \tilde \rho(0,\cdot) \in {\cal D}_n,$ 
then
$$\Vert\rho(t, \cdot) -\tilde\rho(t, \cdot)\Vert_{L^1} \le  
\Vert\rho(0, \cdot) -\tilde\rho(0, \cdot)\Vert_{L^1},$$
for each $0\le t \le T$. 
\end{lemma}

\vspace{5pt}
\n\textsc{Proof.}  Let $\{f_\nu\}$ be a sequence of maps satisfying (${\cal
F}1$), (\ref{appr2}) and (\ref{appr3}) and $\rho_\nu$  be the  admissible
solutions associated to $f_\nu$ such that
\be\rho_\nu (0,\cdot) \equiv \rho (0,\cdot), \qquad \tilde\rho_\nu (0,\cdot)
\equiv \tilde\rho (0,\cdot). \label{in}\ee 
By Lemma \ref{lemma51}, we have 
$\rho_\nu\rightarrow \rho$ and $\tilde\rho_\nu\rightarrow  
\tilde\rho$ in $L^1$.
By Theorem \ref{exstl1} and (\ref{in}), for
each $0\le t \le T$ and $\nu \in\NN$, there results 
$$\Vert\rho(t, \cdot) -\tilde\rho(t, \cdot)\Vert_{L^1} \le
\Vert\rho(t, \cdot) -\rho_\nu(t, \cdot)\Vert_{L^1}
+\Vert\rho_\nu(t, \cdot) -\tilde\rho_\nu(t, \cdot)\Vert_{L^1}+
\Vert\tilde\rho(t, \cdot) -\tilde\rho_\nu(t, \cdot)\Vert_{L^1}\le$$
$$\le\Vert\rho(t, \cdot) -\rho_\nu(t, \cdot)\Vert_{L^1}+\Vert\rho(0, \cdot)
-\tilde\rho(0, \cdot)\Vert_{L^1}+ \Vert\tilde\rho(t, \cdot)
-\tilde\rho_\nu(t, \cdot)\Vert_{L^1} \longrightarrow \Vert\rho(0, \cdot)
-\tilde\rho(0, \cdot)\Vert_{L^1}.$$
So the proof is concluded.\hfill$\Box$

\section[]{Time Dependent Traffic} 
\label{section85}
\indent

In this section we consider a model of traffic including cross lights
and time dependent traffic. The latter means that the choice of drivers
at junctions depends on the period of the day, so during the morning the
traffic flows towards some direction and during afternoon it may change 
towards another direction. This means that the matrix
$A$ depends on time $t$ (see Section \ref{section2}). 

Consider a single junction as in Section \ref{section2} with two incoming
roads and two outgoing ones. Let $\alpha_1=\alpha_1 (t),\> \alpha_2=\alpha_2
(t)$ be two piecewise constant periodic functions such that
\be 
\alpha_{1}(t)\not=\alpha_{2}(t), \qquad \label{matricet}
\ee  
for each $t\ge0.$
Moreover let $\chi_1=\chi_1(t),\chi_2=\chi_2(t)$ be piecewise constant
periodic maps such that
$$\chi_1(t)+\chi_2(t)=1,\qquad\chi_i(t)\in\{0,1\}, \qquad i=1,2$$
for each $t \ge 0$. The two maps represent traffic lights, the value
$0$ corresponding to red light and the value $1$ to green light. 

\begin{definition} \label{defsem} 
Consider $ \rho = \rho (t, x_1,..., x_{4}) = \big(\rho_1 (t, x_1),...,
\rho_{4} (t,x_{4}) \big)$ with bounded variation.
We say that $\rho$ is a solution at the junction $J$ if 
it satisfies (i), (ii), (iv)
of Definition \ref{definitionsoluz} and the following property holds: %
\begin{itemize} %
\item[(v)]$ f(\rho_3 (t, a_3+)) = \alpha_{1}(t)\chi_1(t)
f(\rho_1 (t, b_1-))+\alpha_{2}(t)\chi_2(t)
f(\rho_2 (t, b_2-))$ and \\
$f(\rho_4 (t, a_4+)) =
\big(1-\alpha_{1}(t)\big)\chi_1(t) f(\rho_1 (t,
b_1+))+\big(1-\alpha_{2}(t)\big)\chi_2(t) f(\rho_2 (t, b_2+))$ for each 
$t>0$.
\end{itemize}

\end{definition}

Assume that at time $\bar t$ one of the maps
$\alpha_1(\cdot),\>\alpha_2(\cdot),\>\chi_1(\cdot),\>\chi_2(\cdot)$
jumps, then we have to solve a new Riemann Problem in the junction hence four
waves are generated and
\be N\big(\rho(\bar t+, \cdot)\big) \le
N\big(\rho(\bar t-, \cdot)\big),\label{semeff1}\ee  
\be\tv\big(f(\rho(\bar t+,
\cdot))\big) \le \tv\big(f(\rho(\bar t-, \cdot))\big) + 4 f(\sigma).
\label{semeff2}\ee
Then the map $N\big(\rho(t,\cdot)\big)$ is still non increasing while
$$\tv\big(f(\rho(t_2, \cdot))\big)\le \tv\big(f(\rho( t_1, \cdot))\big) + 4
f(\sigma) \Phi(t_1,t_2),$$
for each $0< t_1\le t_2,$ where
$$ \Phi (t_1,t_2)\doteq \sum_{i=1,2} \Big(\# \big\{t_l\vert t_1<t_l\le t_2, \>
 \chi_i\>\> {\rm jumps \> in}\>\> t_l \big\} +  \# \big\{t_l\vert
t_1<t_l\le t_2, \>  \alpha_i\>\> {\rm jumps \> in}\>\> t_l \big\}\Big).
$$

Therefore, for fixed $T>0$, we have uniform bounds of 
the total variation on the interval $[0,T]$,
and using arguments as in the previous sections we
obtain existence and  stability for the Cauchy Problem.
However, the total variation of $f(\rho)$ does not
depend continuously  on the total variation of the maps
$\alpha_1(\cdot),\>\alpha_2(\cdot)$.  Indeed consider a single junction with
two incoming roads and two outgoing ones  without traffic lights, i.e.
$\chi_i\equiv 1$, and let 
$$\alpha_1 (t)=\left\{\begin{array}{ll} \beta_1 & \textrm{if $0\le t \le \bar
t$},\\  \beta_2 & \textrm{if $\bar t \le t \le T$},
\end{array} \right. , \quad
\alpha_2 (t)=\left\{\begin{array}{ll}
\beta_2 & \textrm{if $0\le t \le \bar t$},\\
\beta_1 & \textrm{if $\bar t \le t \le T$},
\end{array} \right. ,$$
where $0< \beta_2 < \beta_1<{1\over 2}$ and $0 <\bar t < T$.
Consider the initial data $(\rho_{1,0}, \>\rho_{2,0}, \>\rho_{3,0},
\>\rho_{4,0}),$ where
$$f(\rho_{1,0})=f(\rho_{4,0})= f(\sigma),\qquad
f(\rho_{2,0})=f(\rho_{3,0})={{\beta_1}\over {1-\beta_2}}f(\sigma).$$
This is an equilibrium configuration for the choice $\alpha_i=\beta_i$, 
$i=1,2$, hence the solution of the Riemann
Problem is identically equal to the initial data for $0\le t \le \bar t$. At
time $t=\bar t$ we have to solve a new Riemann Problem. Let  $(\hat\rho_{1},
\>\hat\rho_{2}, \>\hat\rho_{3}, \>\hat\rho_{4})$ its solution, there results
(see figure 2)
 $$f(\hat\rho_{2})=f(\hat\rho_{4})= f(\sigma),\qquad
f(\hat\rho_{1})=f(\hat\rho_{3})={{\beta_2}\over {1-\beta_1}}f(\sigma).$$

\begin{center}
\leavevmode \epsfxsize=4in \epsfbox{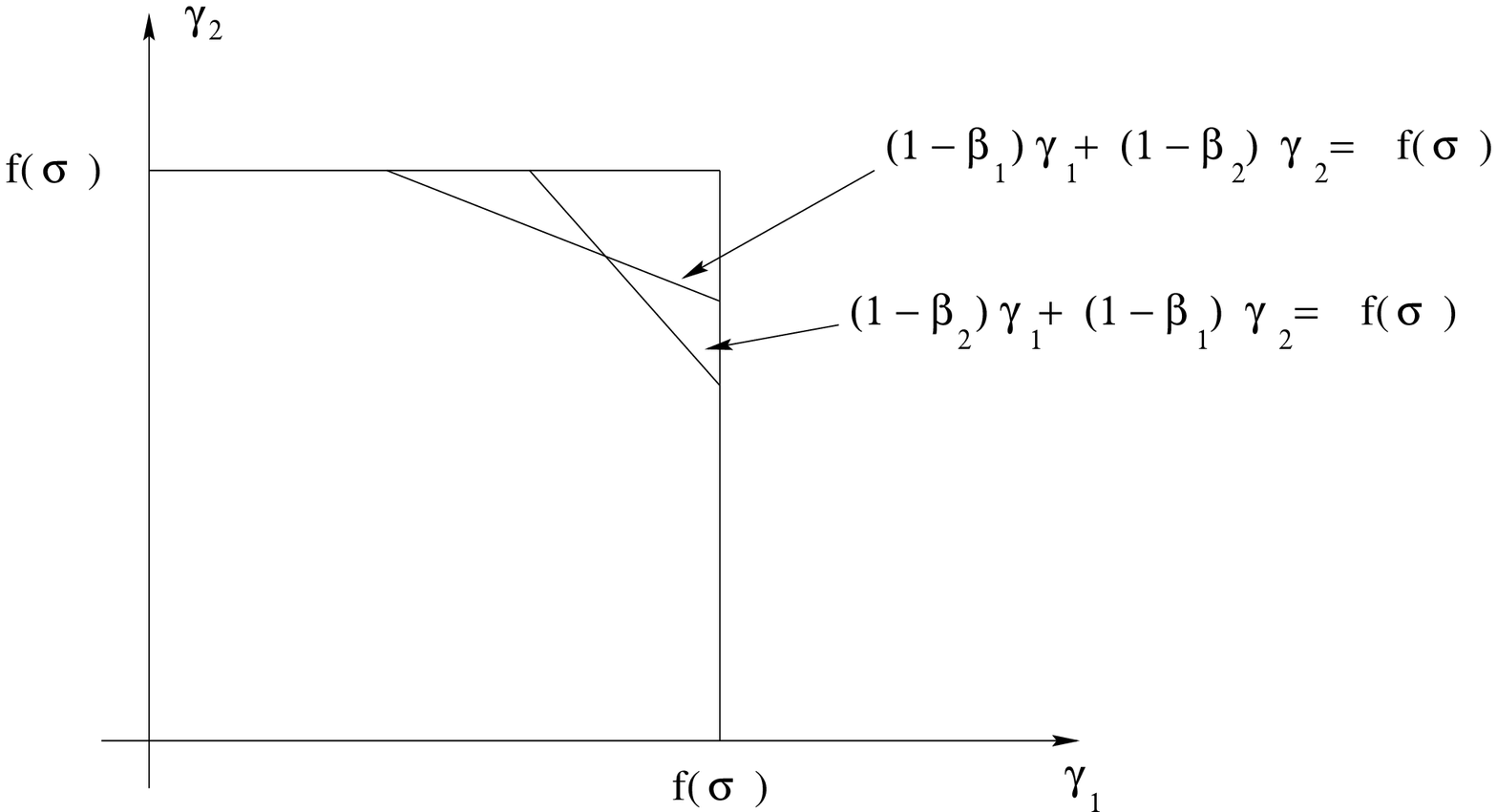}
\end{center}
\begin{center} 
figure 2
\end{center}

\noindent Now, if $\beta_1 \rightarrow  \beta_2, $
then 
$$
\tv\big(\alpha_1 ;[0,T]\big)\longrightarrow 0 , \qquad \tv\big(\alpha_2
;[0,T]\big)\longrightarrow 0,
$$
but
$$\big(f(\rho_{1,0}),\> f(\rho_{2,0})\big)\longrightarrow
\bigg(f(\sigma),\> {{\beta_2}\over {1-\beta_2}}f(\sigma)\bigg),\qquad
\big(f(\hat\rho_{1}),\>f(\hat\rho_{2})\big)\longrightarrow
\bigg({{\beta_2}\over {1-\beta_2}}f(\sigma),\> f(\sigma)\bigg),$$
hence $\tv(f(\rho) ;[0,T])$ is bounded away from zero.

\appendix
\section[]{Appendix: Total Variation of the Fluxes} 
\label{section91}
\indent

In this section we show an example in which the total variation of the
flux increases due to interactions of waves with junctions. 

Consider a single junction with three incoming roads and three outgoing ones,
the matrix
\be 
A\doteq \left(\matrix{ {1\over 2}& {1\over 2}& {1\over 3}\cr
                           {1\over 3}& {1\over 2}& {1\over 2}\cr
                           {1\over 6}&  0        & {1\over6}\cr}
\right)\label{3x3}
\ee
\begin{center}
\leavevmode \epsfxsize=2in \epsfbox{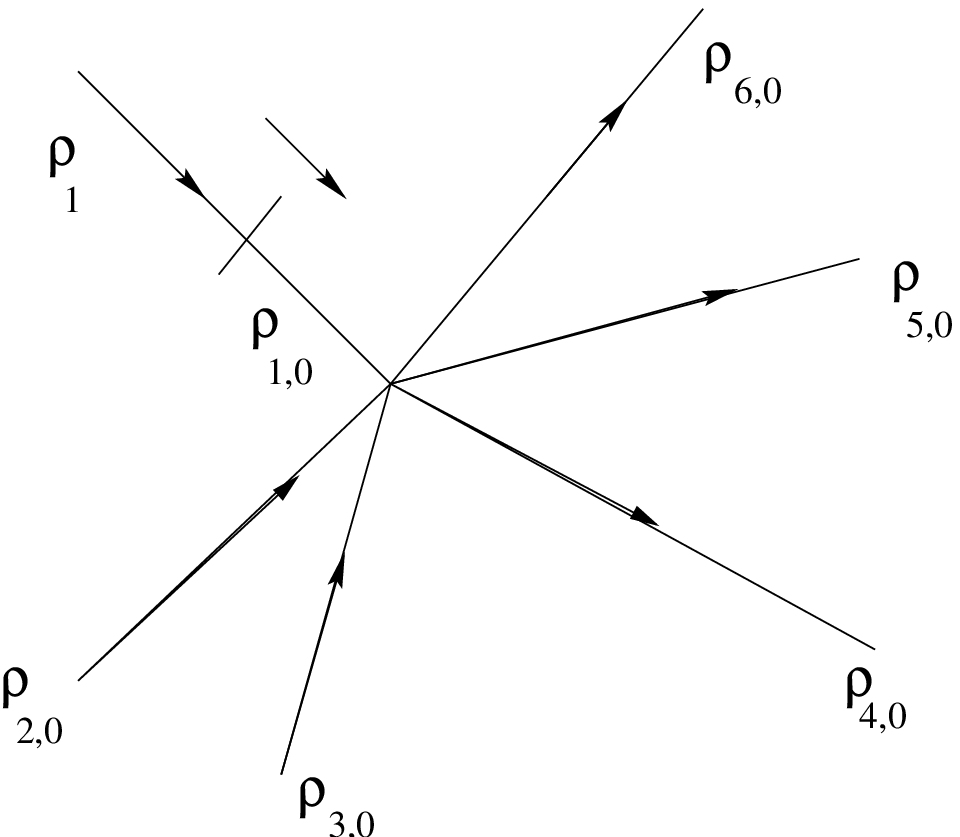}
\end{center}
\begin{center} 
figure 3
\end{center}
and constants $\rho_{1},\>\rho_{1,0},...,\>\rho_{6,0}\in [0,1]$ such that
$$ \rho_{1,0} =\rho_{3,0}=\rho_{4,0}=\rho_{5,0}=\sigma,\>\>
\sigma<\rho_{2,0}<1,\>\> 0<\rho_{6,0},\>\rho_{1}<\sigma,\>\>
f(\rho_{2,0})={1\over 3},\>\> f(\rho_{6,0})={1\over 3}.$$
Assume that $f(\sigma)=1$, then $(\rho_{1,0},...,\>\rho_{6,0})$ 
is an equilibrium configuration and $\rho$ given by 
$$
\rho_{1} (0,x)=
\left\{\begin{array}{ll}  \rho_{1,0}& \textrm{if $x_1\le x \le b_{1} $},\\
 \rho_{1} & \textrm{if $x< x_1$}, \end{array} \right. \quad \rho_i(0,\cdot)
\equiv \rho_{i,0}, \quad i=2,...,6,
$$
is a solution. Moreover the plane 
$${1\over 6}\gamma_1+{1\over 6}\gamma_3=1$$
does not intersect the cube $[0,1]^3$ and the point
$\big(f(\rho_{1,0}),...,\>f(\rho_{6,0})\big)$ is on the intersection of the
planes
$$ {1\over 2}\gamma_1+ {1\over 2}\gamma_2+ {1\over 3} \gamma_3=1, \qquad
{1\over 3}\gamma_1+ {1\over 2}\gamma_2+ {1\over 2}\gamma_3=1,$$
that is the line described by the map
\be 
\gamma_1 \mapsto \Big(\gamma_1,2- {5\over 3}\gamma_1, \gamma_1\Big).
\label{retta}
\ee

At some time say $\bar t$ the wave $(\rho_1,\rho_{1,0})$ interacts with the
junction.  Let $(\hat\rho_{1}, ....,\hat\rho_{6})$ be the solution of the
Riemann Problem at the junction 
for the data $(\rho_{1},\>\rho_{2,0},...,\>\rho_{6,0})$.
Since the map $E$ increases on the line described by (\ref{retta}), the
point $\big(f(\hat\rho_{1}), ....,f(\hat\rho_{6})\big)$ is on the curve 
(\ref{retta}) and
$$
f(\hat\rho_{1})=f(\hat\rho_{3})=f(\rho_1),\>\> f(\hat\rho_{2})= 2-{5\over
3}f(\rho_1),\>\>f(\hat\rho_{4})=f(\hat\rho_{5})=f(\sigma),\>\>
f(\hat\rho_{6})={1\over 3}f(\rho_1).
$$ 
We get
$$
\tv\big(f(\rho(\bar t -, \cdot))\big)=
1-f(\rho_1),
$$
while 
$$\tv\big(f(\rho(\bar t +, \cdot))\big)=
4\big(1-f(\rho_1)\big)> \tv\big(f(\rho(\bar t -, \cdot))\big).
$$

\section[]{Appendix: Total Variation of the Densities} 
\label{section9}
\indent

Consider a junction $J$ with two incoming roads and two outgoing ones that we
parameterize with the intervals $]-\infty,\> b_1],\>]-\infty,\> b_2],\> [a_3,\>
+\infty [,\>[a_4,\>+\infty [$ respectively. 
Fix the constants $\alpha_1,\> \alpha_2$ such that
$0<\alpha_1<\alpha_2<1/2$ and set
$$
\alpha_{3,1}=\alpha_1, \quad \alpha_{3,2}=\alpha_2, \quad
\alpha_{4,1}=1-\alpha_1, \quad\alpha_{3,2}=1-\alpha_2.
$$

\begin{center}
\leavevmode \epsfxsize=2in \epsfbox{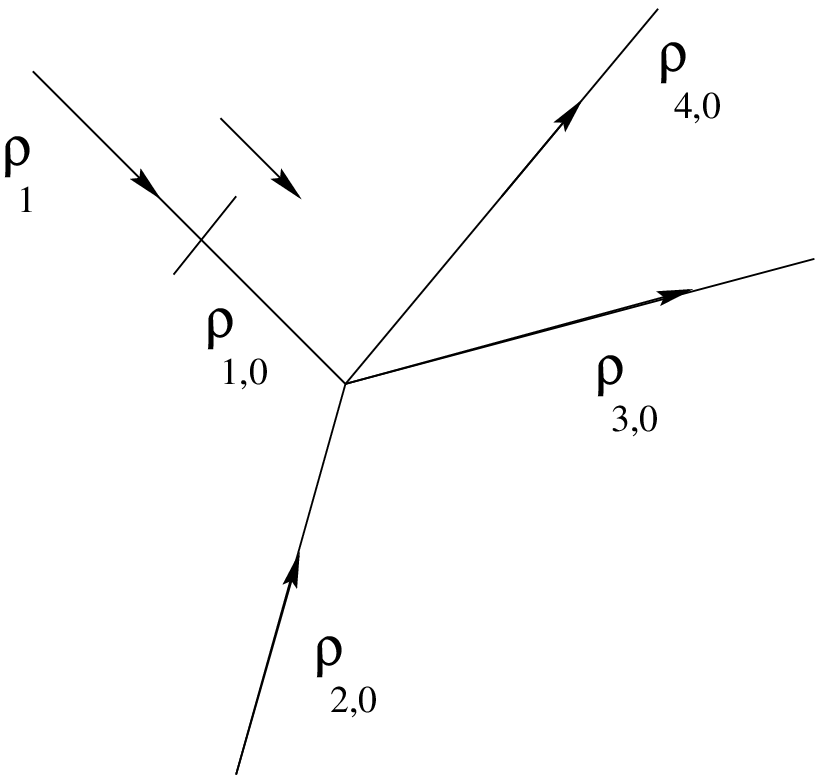}
\end{center}
\begin{center} 
figure 4
\end{center}

\n Define a solution $\rho$ by
\be
\rho_{1} (0,x)= \left\{\begin{array}{ll}
\rho_{1,0}& \textrm{if $x_1\le x \le b_{1} $},\\
\rho_{1} & \textrm{if $x< x_1$},
\end{array} \right.  \quad \rho_{2} (0,x)=\rho_{2,0},\quad \rho_{3}
(0,x)=\rho_{3,0}, \quad \rho_{4}(0,x)=\rho_{4,0},\label{datoiniziale}\ee
where $\rho_{1},\>\rho_{1,0},\>\rho_{2,0},\>\rho_{3,0},\>\rho_{4,0}$ are
constants such that 
\be
\sigma \le \>\rho_{2,0}, \>\rho_{3,0}\le 1, \quad 0\le
\rho_1\le\sigma,\quad \rho_{1,0}= \rho_{4,0}=\sigma,\label{dati}
\ee
$$
f(\rho_{1,0})=f(\rho_{4,0})= f(\sigma),\quad
f(\rho_{2,0})=f(\rho_{3,0})={{\alpha_1}\over {1-\alpha_2}}f(\sigma),
$$
so $(\rho_{1,0}, \>\rho_{2,0}, \>\rho_{3,0}, \>\rho_{4,0})$ is an equilibrium
configuration. 

After some time the wave $(\rho_1,\rho_{1,0})$ interacts with the junction.
Let $(\hat\rho_{1}, \>\hat\rho_{2}, \>\hat\rho_{3}, \>\hat\rho_{4})$
be the solution of the Riemann Problem in the junction 
for the data $(\rho_{1}, \>\rho_{2,0},
\>\rho_{3,0}, \>\rho_{4,0})$. By (\ref{datoiniziale}) and
(\ref{dati}) there results
$$
f(\hat\rho_{1})=f(\rho_{1}), \quad f(\hat\rho_{2})= {{f(\sigma)-
(1-\alpha_1) f(\rho_{1})}\over {1-\alpha_2}},
$$
$$
f(\hat\rho_{3})= {{\alpha_1 -\alpha_2}\over {1-\alpha_2}}f(\rho_{1}) +
{{\alpha_2}\over {1-\alpha_2}}f(\sigma), \quad f(\hat\rho_{4})=f(\sigma)
$$
and 
\be 0\le \hat\rho_3 \le\sigma\le \hat\rho_2 \le 1.
\label{datisb}\ee
\vspace{2pt}

\begin{center}
\leavevmode \epsfxsize=4in \epsfbox{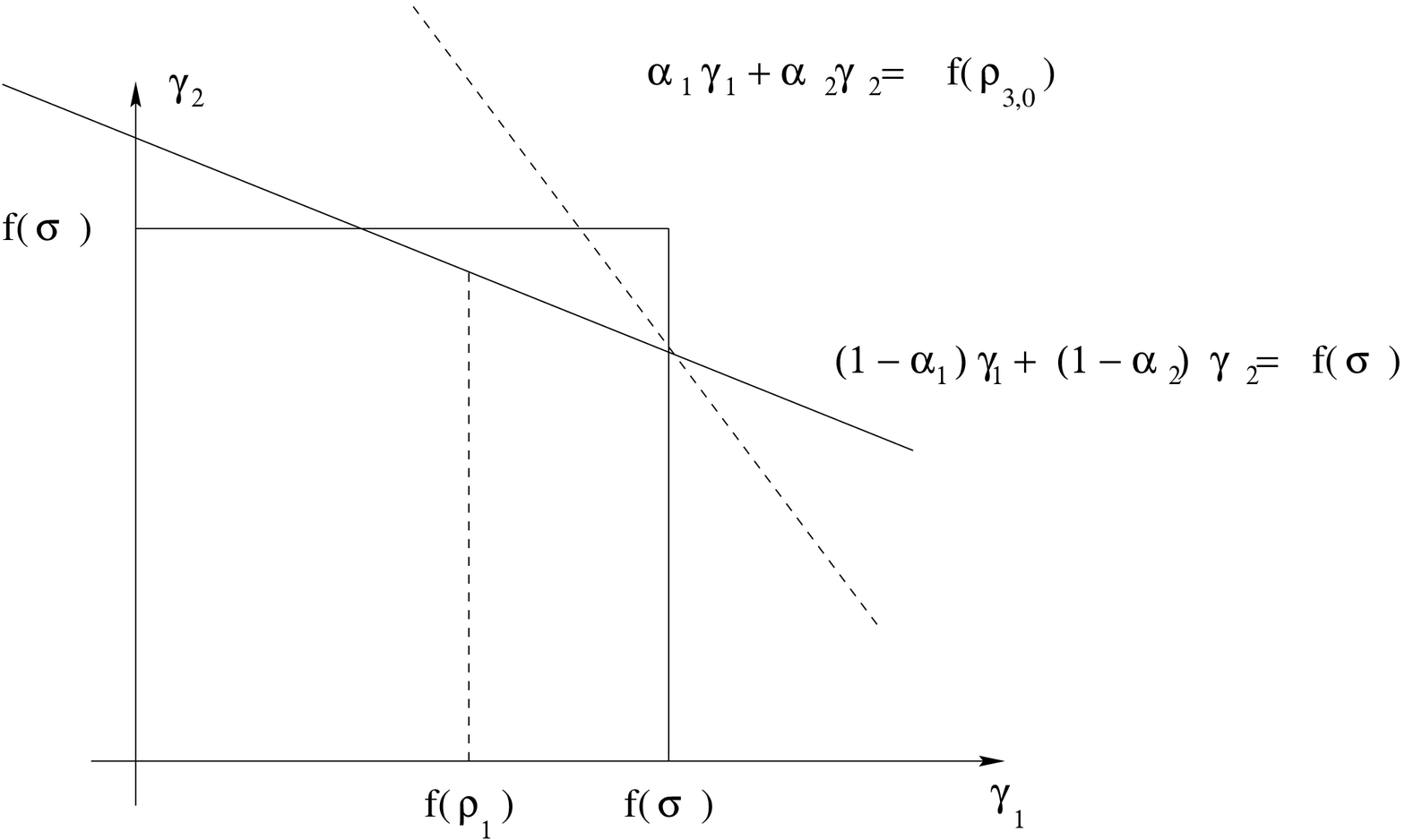}
\end{center}
\begin{center}
figure 5
\end{center} 
Therefore if $\rho_1\rightarrow \rho_{1,0} =\sigma$ then
$$f(\hat\rho_3)\longrightarrow{{\alpha_1}\over
{1-\alpha_2}}f(\sigma)=f(\rho_{3,0}),$$
and by (\ref{datisb}), (\ref{dati}) we have
$\hat\rho_3\rightarrow \tau(\rho_{3,0})$.
Therefore we are able to create on the third road a wave with strength
bounded away from zero using an arbitrarily small wave on the first one.

\vsp
\centerline{\Ack Acknowledgements}
\vsp

\noindent
{\ack 
The authors
would like to thank Prof. Rinaldo M. Colombo 
for useful discussions. 
}

%
%
\newcommand{\auth}{\textsc}
\newcommand{\tit}{\textrm}
\newcommand{\jou}{\textit}

\end{document}